\documentclass[a4paper,11pt]{amsart}
\usepackage{amsthm,amsfonts,amsbsy,amssymb,amsmath,amscd}
\usepackage[all]{xy}
\usepackage[colorlinks=true]{hyperref}
\usepackage{pgf, tikz, tikz-cd}

\DeclareMathOperator{\vol}{vol}

\newtheorem{teo}{Theorem}

\theoremstyle{definition}
\newtheorem{df}[teo]{Definition}

\theoremstyle{remark}

\numberwithin{equation}{section}

\begin{document}

\title[On Canonical threefolds near the Noether line]{On Canonical threefolds\\ near the Noether line}

\author{Roberto Pignatelli}
\date{\today}
	
\address{Dipartimento di Matematica dell'Universit\`a di Trento, via Sommarive 14, I-38123, Trento (Italy)}
\email{Roberto.Pignatelli@unitn.it}


\thanks{This note originates from a series of seminars I conducted in 2023 and 2024 in Bedlewo (Poland), Pipa (Brazil), Taipei (Taiwan), Shanghai (China), Pavia, and Pisa (Italy). It has greatly benefited from the numerous contributions of the knowledgeable audiences I encountered on each of these occasions. I warmly thank all the participants as well as the colleagues who honored me with their invitations.
The author was partially supported by the "National Group for Algebraic and Geometric Structures, and their Applications" (GNSAGA - INdAM) and by the European Union - Next Generation EU, Mission 4 Component 2 - CUP E53D23005400001}

\maketitle

\setcounter{tocdepth}{1}

Let $X$ be a canonical $n$-fold, that is a complex irreducible projective variety of dimension $n$ whose singularities are at worst canonical and whose canonical class $K_X$ is ample. If $n=2$ we will use the word "surface" instead of the word $2$-fold. 

As usual we denote by $p_g:=h^0(X,K_X)$ its geometric genus and by $K_X^n$ its canonical volume.
The volume of a canonical surface is an integer respecting the M. Noether inequality $K^2 \geq 2p_g-4$.  

The moduli space of the canonical surfaces $X$ with $K_X^2=2p_g(X)-4$ has been described by Horikawa.
According to Horikawa's description, if $K_X^2=2p_g(X)-4$ and $p_g \geq 7$, then there is a fibration $f\colon X \rightarrow {\mathbb P}^1$ with fibres of genus $2$.
Let $f \colon X \rightarrow B$ be a genus $2$ fibration, that is a surjective morphism onto a smooth projective curve $B$ whose general fibre is a curve of genus $2$. The canonical ring of a fibre $F$, the ring $\bigoplus_d H^0(F,dK_F)$, is of one of the following two forms:
\begin{align*}
\text{\it $2$-connected fibres}&&&&\text{\it $2$-disconnected fibres \ \ \ }\\
\frac{{\mathbb C}[x_0,x_1,z]}{f_6(x_0,x_1,z)}&&
\text{ or }&&
\frac{{\mathbb C}[x_0,x_1,y,z]}{f_2(x_0,x_1),f_6(x_0,x_1,y,z)}
\end{align*}
with $\deg x_j=1$, 
$\deg y=2$, 
$\deg z=3$, 
$\deg f_d=d$.

Moreover (see \cite{CPlowgenus})
\[
K^2_X=2p_g\left( X \right)-4+6b-2h^1(X,{\mathcal O}_X)+\deg \tau
\]
where $b$ is the genus of the base curve $B$, and $\tau$ is an effective divisor on $B$ supported on the image of the $2$-disconnected  fibres.
By an inequality of Jongmans and Debarre, if $K^2_X < 2p_g$, then $X$ is regular (that means $h^1(X,{\mathcal O}_X)=0$), which in turn implies $b=0$.
In particular  
\[
K^2_X=2p_g\left( X \right)-4 \Leftrightarrow h^1(X,{\mathcal O}_X)=0 \text{ and } \tau=0.
\]

We can then state the following.
\begin{df}
A simple genus $2$ fibration is a projective variety $X$ with a morphism $f$ on a smooth projective curve $B$ such that
\begin{enumerate}
\item  all singularities of $X$ are canonical;
\item $K_X$ is $f$-ample;
\item for {\em all} $p \in B$, the canonical ring of the fibre $X_p=f^{-1}p$ is of the form
$\frac{{\mathbb C}[x_0,x_1,z]}{f_{6}(x_0,x_1,z)}$, with $\deg x_i=1$, $\deg z=3$, $\deg f_{6}=6$
\end{enumerate}
\end{df}
and the mentioned results give the following
\begin{teo}\label{thm: Horikawa fibrations} 
Let $X$ be a projective variety of dimension $2$ with at worst canonical singularities and $p_g\geq 7$.
Then $X$ is a regular simple genus $2$ fibration if and only if $X$ is a canonical surface and $K^2_X=2p_g-4$. 
\end{teo}

If $n=3$ \cite{CCJ,CCJAdd} proved that  the {\it Noether inequality for canonical $3$-folds} 
\[
K_X^3\geq \frac43 p_g(X)- \frac{10}3
\]
holds, unless (possibly)  $5\leq p_g \leq 10$. 
The first step of the proof is showing that, if $K^3_X$ is not much bigger than 
$ \frac43 p_g(X)- \frac{10}3$, then there exists a fibration 
$ f \colon X \dashrightarrow {\mathbb P}^1$ whose general fibre is a canonical surface with canonical volume $1$ and geometric genus $2$. These are classically known to be the hypersurfaces of degree $10$ in ${\mathbb P}(1,1,2,5)$ with at worst canonical singularities.

This, together with Theorem \ref{thm: Horikawa fibrations}, inspired the following definition
\begin{df}\cite{CP}
A simple fibration in $(1,2)$-surfaces is a projective variety  with a morphism $f$ on a smooth projective curve $B$ such that
\begin{enumerate}
\item  all singularities of $X$ are canonical;
\item $K_X$ is $f$-ample;
\item for {\em all} $p \in B$, the canonical ring of the fibre $X_p=f^{-1}p$ is of the form $\frac{{\mathbb C}[x_0,x_1,y,z]}{f_{10}(x_0,x_1,y,z)}$,
with $\deg x_i=1$, $\deg y=2$, $\deg z=5$, $\deg f_{10}=10$.
\end{enumerate}
\end{df}

In \cite{CHPZ1} we generalize Theorem \ref{thm: Horikawa fibrations} to dimension $3$ as follows.
\begin{teo}\label{thm: main}
Let $X$ be a projective variety of dimension $3$ with at worst canonical singularities and $p_g\geq 23$.
Then $X$ is a Gorenstein regular simple fibration in $(1,2)$-surfaces if and only if $X$ is a canonical $3$-fold and $K^3_X=\frac 43 p_g-\frac{10}3$. 
\end{teo}
One implication was proven in \cite{CP} where the other implication was conjectured. 
The inequality $p_g \geq 23$ is sharp, since there are simple fibrations in $(1,2)$-surfaces with $p_g=13,16,19,22$ and $K_X$  not ample: the map from these threefolds onto their canonical model is a small contraction, contracting a rational curve contained in the singular locus of $X$.

This allowed a complete classification of the moduli spaces ${\mathcal M}_{p_g}$ of the canonical threefolds with $p_g \geq 11$ and $K_X^3\geq \frac43 p_g(X)- \frac{10}3$. It was already known that ${\mathcal M}_{p_g}$ is empty unless $p_g+2$ is divisible by $3$.
In \cite{CHPZ1} we prove 
\begin{teo}
If $p_g \geq 11$ and $p_g+2$ is divisible by $3$ then ${\mathcal M}_{p_g}$ decomposes as union of $\left\lfloor \frac{p_g+6}4\right\rfloor$ unirational strata. 
The number of irreducible components of ${\mathcal M}_{p_g}$ is at most $\left\lfloor \frac{p_g+6}4\right\rfloor$ and at least  $\left\lfloor \frac{p_g+6}4\right\rfloor-\left\lfloor \frac{p_g+8}{78}\right\rfloor$.
 \end{teo}
 
 It is worth mentioning that, for all $p_g$ only one or two strata contain smooth $3$-folds. The remaining strata provide then examples of $3$-folds with smoothable singularities that are not globally smoothable.
 
The following results 
\begin{teo}[\cite{CP}]\label{thm: simple are 2-Gorenstein}
Let $f \colon X \rightarrow B$ be a simple fibration in $(1,2)$-surfaces. Then $X$ is $2$-Gorenstein.
If $X$ is regular then $K^3_X=\frac43 p_g - \frac{10}3 +\frac{N}6$, $N \in {\mathbb N}$. Moreover $X$ is Gorenstein if and only if
 $N=0$.
\end{teo}

\begin{teo}[\cite{HZ}]\label{thm: three lines}
If $X$ is a canonical $3$-fold with $p_g(X) \geq 11$ and $K^3_X \leq \frac43 p_g \left( X \right) -\frac{10}3 +\frac{2}6$ then $X$ is 2-Gorenstein and $K^3_X = \frac43 p_g \left( X \right) -\frac{10}3 +\frac{N}6$ with $N=0,1,2$. Moreover $X$ is Gorenstein if and only if $N=0$.
\end{teo}

lead to the natural

{\bf Question:} {\em 
Find the maximal $\epsilon >0$ (if it exists) such that if $p_g>>1$ then all threefolds $X$ with 
$\vol(X)< \frac43 p_g \left( X \right)-\frac{10}3 + \epsilon $ are regular simple fibrations in $(1,2)-$surfaces.}

Theorems \ref{thm: main} and \ref{thm: three lines} show $\epsilon \geq \frac16$. Proving  $\epsilon \geq \frac26$ would lead to a complete classification of the moduli space of canonical threefolds {\em on the second Noether line}. 

On the opposite side we can prove that
$
\epsilon \leq \frac46
$
by constructing  (\cite{CHPZ2}), for each positive integer $p_g$ divisible by $3$, canonical $3$-folds with $K^3=\frac43p_g+\frac{10}3+\frac46$ of Gorenstein index $3$, which implies by Theorem \ref{thm: simple are 2-Gorenstein} that they are not simple fibrations in $(1,2)$-surfaces.


\bibliography{Ref_firstNoetherline}
\bibliographystyle{amsalpha}

\end{document}